\documentclass[final,3p]{elsarticle} 
\usepackage{eurosym}
\usepackage{amsfonts}
\usepackage{amsmath}
\usepackage{amssymb}
\usepackage{graphicx}
\usepackage[tableposition=top]{caption}
\usepackage{subcaption}
\usepackage{float}
\usepackage{setspace}
\usepackage{placeins}
\usepackage{mathrsfs}

\newsavebox \foobox
\newlength{\foodim}

\graphicspath{{C:/Figures/}}

\usepackage{multirow}

\makeatletter
\def\ps@pprintTitle{%
 \let\@oddhead\@empty
 \let\@evenhead\@empty
 \let\@oddfoot\@empty
 \let\@evenfoot\@empty
}
\makeatother

\begin{document}
\begin{frontmatter}
\title{A computationally efficient fractional predictor–corrector approach involving the Mittag–Leffler kernel
}
\author{Sami Aljhani$^{a}$}
\address{$(a)$Department of Mathematics \& Statistics, College of Science in Yanbu , Taibah University,Yanbu Govemorate, Saudi Arabia.\\ 
	(*) Corresponding author, email: sjhani@taibahu.edu.sa}
\begin{abstract}
In this paper, based on Newton interpolation we have proposed a numerical scheme of predictor-corrector type in order to solve fractional differential equations with the fractional derivative involving the Mittag-Leffler function. We have added an auxiliary midpoint in each sub-interval, this allows us to use a piecewise quadratic Newton interpolation to derive the corrector scheme. The derivation of the schemes for the midpoint and the predictor is done by means of a piecewise linear Newton interpolation. We present some illustrative examples for initial value problems that involve fractional derivatives in the sense of Atangana-Baleanu. The results of numerical experiments show that the proposed scheme is a powerful technique to handle fractional differential equations with nonlinear terms that involve operators of Atangana-Baleanu type. Moreover, the proposed method significantly improves the numerical accuracy in comparison with other methods.
\end{abstract}
\end{frontmatter}

\section{Introduction}
Although fractional calculus has a long history in mathematics, researchers have only recently turned their attention to its practical applications, owing to the ability of fractional derivatives to model phenomena with memory and hereditary effects \cite{alquran2018novel,douaifia2020asymptotic,Osman2017}.
Several distinct formulations of fractional derivatives and integrals have been introduced, most notably those of Riemann–Liouville, Caputo, and Grünwald–Letnikov  \cite{magin2006fractional}.
Later, Caputo and Fabrizio \cite{caputo2015new} proposed a new definition of the fractional derivative involving an exponential decay kernel, which inspired Atangana and Baleanu \cite{atangana2016new}  to introduce a more general formulation characterized by a nonsingular, nonlocal kernel based on the Mittag–Leffler function.

Usually, it is hard to find exact and direct solutions in fractional calculus, so we rely on numerical and approximate methods because they give results that are very close to the real ones. In recent years, there’s been a lot of attention on these numerical methods since they help solve fractional differential equations that can not easily be handled by hand  \cite{douaifia2021newton,Osman2017}.

Since these methods are used in many areas including business, arts, medicine, life, and social sciences and computers are an essential part of them, they often provide more accurate results than standard analytical solutions. Moreover, numerical analysis helps us understand how real-world models behave. That is why developing and applying accurate numerical techniques has become so important for solving real-life problems \cite{atangana2016new, Toufik2017, ullah2020modeling,Osman2017}.

A long time ago ,more than two thousand years ago,  people started using simple numerical methods. They began with an easy idea called linear interpolation. Later on, scientists created more advanced ways to solve similar problems, like Newton’s method, Euler’s method, Gaussian elimination, and Lagrange interpolation \cite{werner1984polynomial, dimitrov1994note, yang2016visualizing, krogh1970efficient}.

In the past few years \cite{diethelm1998fracpece, diethelm2004detailed, djida2017numerical, baleanu2018nonlinear}, many researchers have focused on creating numerical methods to solve fractional models that come from ordinary differential equations. One well-known method is the Adams–Bashforth method, which is popular because it uses a simple idea called Lagrange interpolation \cite{jain2018numerical, zhang2018decoupled}. Recently, another version of this method was introduced that uses Newton’s quadratic interpolation instead \cite{ douaifia2021newton,atangana2021new}.

For fractional differential equations with the Caputo fractional derivative, Diethelm et al.  \cite{diethelm2002predictor} used the definition of the Caputo derivative to produce a Predictor–Corrector approach for a fractional differential equation.  To perform this method,  convert the fractional ordinary differential equation into a Volterra Integral Equation,  then use the Rectangle rule for the Predictor and the Trapezoid rule for the Corrector steps.  The convergence order of the predictor-corrector approach was $min(2, 1 + \alpha)$, where $\alpha > 1$.  They also talked about how to improve accuracy and reduce processing costs by using techniques like Richardson extrapolation and the short memory principle.  The author showed that the smoothness properties of the given function $f$ and the unknown solution $y$ play an important role.  The researchers find that the smoothness of one of these functions will imply non-smoothness of the other unless some special conditions are fulfilled.  When compared to Diethelm et al approach’s, Deng and Li \cite{nguyen2017high} described how to use the short memory principle combined with the predictor–corrector approach to reduce the computational cost.  They provided an improved method of the Predictor–Corrector algorithm in Diethelm et al. with an accuracy increased to $min(2, 1 + 2 \alpha)$,  where $\alpha > 1$,  and half the computational cost.  The authors present two examples to illustrate the performance of the proposed predictor-corrector schemes. Thereafter,  researcher foucos proposes unique techniques that will result in better convergence.  For example,  In \cite{li2011numerical},  the authors looked at an algorithm that uses the Simpson-formats to calculate the intermediate nodes in order to improve the accuracy and computational complexity of a specific type of predictor algorithm.

For fractional differential equations that use the Atangana–Baleanu derivative, Toufik and Atangana \cite{Toufik2017} used a simple two-step method based on Lagrange interpolation to find the solution, and it worked good for many real-life problems.
Later, Baleanu and his team \cite{baleanu2018nonlinear}  introduced another method called the Adams predictor–corrector. They studied if the equations they worked on really have a solution and if that solution is unique when using the new AB-Caputo form. They also tried two ways to solve these equations and checked how accurate and steady they are using something called Grönwall’s inequality.
Their results showed that the predictor–corrector method gives better and faster results than the Euler method when using the same step size. They also noticed that the Euler method gives a first-level accuracy, while the predictor–corrector method gives a higher, second-level accuracy. In general, both methods are quite good and stable for solving this type of fractional equations.

Many other numerical methods have also been created to solve fractional differential equations that use the Atangana–Baleanu derivative, and you can find several examples of them in the research works. \cite{sadeghi2020operational,ganji2020new,abdeljawad2020analysis,Aljhani2025}
.
So, the goal of this study is to show an easy and effective numerical method to solve fractional problems that include the Atangana–Baleanu fractional derivative.
 \cite{atangana2016new}.

This paper is organized as follows. In Section 2, we explain the main ideas and symbols, like the Atangana–Baleanu derivative and its related integral. In Section 3, we show a new and easy method to solve problems that involve fractional derivatives. In Section 4, we give examples to show how the method works and compare it with some other methods \cite{Toufik2017,baleanu2018nonlinear}.

\section{\textbf{Preliminaries and Notations}}

In this section, we give a short overview of some basic ideas and operators in fractional calculus that are needed in this study.  

The fractional derivative of Atangana--Baleanu in the Caputo sense is written as~\cite{atangana2016new,Aljhani2020numerical}:
\[
{}^{ABC}_{0}D^{\alpha}_{t} f(t)
= \frac{M(\alpha)}{1-\alpha}\!\int_{0}^{t} f'(s)\,
E_{\alpha}\!\left(\frac{-\alpha}{1-\alpha}(t-s)^{\alpha}\right)\!ds,
\qquad 0<\alpha\le1.
\]

Here, $\alpha$ denotes the fractional order and $M(\alpha)>0$ is a normalization function that satisfies $M(0)=M(1)=1$.  
The term $E_{\alpha}(\cdot)$ represents the \textit{Mittag--Leffler function} of order $\alpha$, given by
\[
E_{\alpha}(z)=\sum_{k=0}^{\infty}\frac{z^{k}}{\Gamma(\alpha k+1)},\qquad \Re(\alpha)>0.
\]

The well-known \textit{Gamma function} is defined as
\[
\Gamma(\alpha)=\int_{0}^{\infty} t^{\alpha-1}e^{-t}\,dt,\qquad \Re(\alpha)>0.
\]

The fractional integral related to the ABC operator, which has a non-local kernel and does not have singularities at $t=s$,  
can be expressed as~\cite{atangana2016chaos}:
\[
{}^{ABC}_{0}I^{\alpha}_{t} f(t)
= \frac{1-\alpha}{M(\alpha)}f(t)
+\frac{\alpha}{M(\alpha)\Gamma(\alpha)}\int_{0}^{t} f(s)(t-s)^{\alpha-1}ds.
\]

When $\alpha=0$, we obtain the original function $f(t)$,  
and when $\alpha=1$, the expression reduces to the classical integral.

\section{Improved Numerical method for fractional differential involved Mittag-Leffler function }
We consider the initial-value problem with Atangana-Baleanu derivative 
\begin{equation}\label{eq1}
\begin{cases}
&{}^{ABC}_{\quad \, 0}D^{\alpha}_{t} y\left(t\right)= f\left(t,y\left(t\right)\right),   \qquad   \qquad    \qquad       0< t< T < \infty, \\ \\
& \qquad  \quad \, y\left(0\right)= y_{0},   
\end{cases}
\end{equation}
where $f$ is a smooth nonlinear function that guarantees the existence of a unique solution $y\left(t\right)$ for (\ref{eq1})  and the fractional order $\alpha \in \left(0,1\right)$. In order to obtain a predictor-corrector numerical scheme that solves (\ref{eq1}),  we apply the Atangana-Baleanu integral,  we convert the above equation into 
\begin{align}
&y\left(t\right)-y\left(0\right)=\frac{1-\alpha}{AB(\alpha)}f\left(t,y\left(t\right)\right)+\frac{\alpha}{AB(\alpha)\Gamma(\alpha)}\int^{t}_{0}f\left(s,y\left(s\right)\right)\left(t-s\right)^{\alpha-1}ds
\end{align}
at the point  $t_{m+1}=\left(m+1\right)\Delta t$ ,  we have 
\begin{multline}\label{eq2}
y\left(t_{m+1}\right)=y\left(0\right)+ \frac{1-\alpha}{AB(\alpha)}f\left(t_{m+1},y_{m+1}\right)+\frac{\alpha}{AB(\alpha)\Gamma(\alpha)}\int^{t_{m+1}}_{0}f\left(s,y\left(s\right)\right)\left(t_{m+1}-s\right)^{\alpha-1}ds\\
\approx y\left(0\right)+ \frac{1-\alpha}{AB(\alpha)}f\left(t_{m+1},y_{m+1}\right)+\frac{\alpha}{AB(\alpha)\Gamma(\alpha)} \sum^{m}_{j=0} \Bigg\{  \int^{t_{j+1}}_{t_{j}}\tilde{f_{j}}\left(s,y\left(s\right)\right)\left(t_{m+1}-s\right)^{\alpha-1}ds \Bigg\} ,\quad
\end{multline}
where $t_{0}=0$ and $\tilde{f_{j}}\left(s,y\left(s\right)\right)$ is the piecewise quadratic Newton interpolation (three steps) given by 
\begin{multline}
\tilde{f_{j}}\left(s,y\left(s\right)\right)=N_{j}\left(s\right)= f\left(t_{j},y_{j}\right)\\
+\frac{2 \left[f\left(t_{j+\frac{1}{2}},y_{j+\frac{1}{2}}\right)-f\left(t_{j},y_{j}\right)\right]}{\Delta t} \left(s-t_{j}\right)\qquad \qquad \qquad \qquad \\
+\frac{2f\left(t_{j},y_{j}\right)-4f\left(t_{j+\frac{1}{2}},y_{j+\frac{1}{2}}\right)+2f\left(t_{j+1},y_{j+1}\right)}{\left(\Delta t\right)^{2}}\left(s-t_{j}\right)\left(s-t_{j+\frac{1}{2}}\right), \,\,\,\,\,\qquad \quad
\end{multline}
where 
\begin{align*}
&t_{j+\frac{1}{2}}=\left(j+\frac{1}{2}\right)\Delta t .
\end{align*}
Thus replacing this equalities into (\ref{eq2}) lead to the following scheme
\begin{multline}
y\left(t_{m+1}\right)=y_{0}+ \frac{1-\alpha}{AB(\alpha)}f\left(t_{m+1},y_{m+1}\right)+\frac{\alpha}{AB(\alpha)\Gamma(\alpha)} \\
\times \sum^{m}_{j=0} \Bigg[f\left(t_{j},y_{j}\right) \int^{t_{j+1}}_{t_{j}}\left(t_{m+1}-s\right)^{\alpha-1}ds \qquad  \qquad  \qquad  \qquad  \qquad  \qquad  \qquad  \qquad  \qquad   \qquad \quad  \,\,\,\,  \,\,\, \\ 
+\frac{2 \left[f\left(t_{j+\frac{1}{2}},y_{j+\frac{1}{2}}\right)-f\left(t_{j},y_{j}\right)\right]}{\Delta t}\int^{t_{j+1}}_{t_{j}}\left(s-t_{j}\right)\left(t_{m+1}-s\right)^{\alpha-1}ds \qquad \qquad \qquad \qquad \quad \quad  \,\, \,\, \,\,\, \\
+\frac{2f\left(t_{j},y_{j}\right)-4f\left(t_{j+\frac{1}{2}},y_{j+\frac{1}{2}}\right)+2f\left(t_{j+1},y_{j+1}\right)}{\left(\Delta t\right)^{2}}\int^{t_{j+1}}_{t_{j}}\left(s-t_{j}\right)\left(s-t_{j+\frac{1}{2}}\right)\left(t_{m+1}-s\right)^{\alpha-1}ds\Bigg],
\end{multline}
which can be simplified and rearranged to the form 
\begin{multline}\label{eq3}
y_{m+1}=y_{0}+ \frac{1-\alpha}{AB(\alpha)}f\left(t_{m+1},y_{m+1}\right)\\
\quad \quad\,\,+\frac{\alpha}{AB(\alpha)\Gamma(\alpha)} \sum^{m}_{j=0}f\left(t_{j},y_{j}\right) \int^{t_{j+1}}_{t_{j}}\left(t_{m+1}-s\right)^{\alpha-1}ds \qquad \qquad \qquad \qquad \qquad \qquad \quad \qquad  \quad \qquad \quad   \\ 
\quad \quad\,\,+\frac{2 \alpha}{AB(\alpha)\Gamma(\alpha)} \sum^{m}_{j=0}\frac{f\left(t_{j+\frac{1}{2}},y_{j+\frac{1}{2}}\right)-f\left(t_{j},y_{j}\right)}{\Delta t}\int^{t_{j+1}}_{t_{j}}\left(s-t_{j}\right)\left(t_{m+1}-s\right)^{\alpha-1}ds \quad \qquad  \quad\,\,\,\,  \qquad \quad \\ 
\quad \qquad+\frac{\alpha}{AB(\alpha)\Gamma(\alpha)} \sum^{m}_{j=0} \frac{\left(2f\left(t_{j},y_{j}\right)-4f\left(t_{j+\frac{1}{2}},y_{j+\frac{1}{2}}\right)+2f\left(t_{j+1},y_{j+1}\right)\right)}{\left(\Delta t\right)^{2}} \qquad   \qquad \qquad \quad \qquad  \qquad  \quad \,\, \\
\times  \int^{t_{j+1}}_{t_{j}}\left(s-t_{j}\right)\left(s-t_{j+\frac{1}{2}}\right)\left(t_{m+1}-s\right)^{\alpha-1}ds. \qquad \qquad \qquad  \qquad \qquad \qquad \qquad \qquad \quad \quad \quad \,\,\,\,\,\,
\end{multline}
The integral in (\ref{eq3}) can be calculated as 
\begin{align}
&\int^{t_{j+1}}_{t_{j}} \left( t_{m+1}-s\right)^{\alpha-1}ds=
\frac{\left(\Delta t\right)^{\alpha}}{\alpha} \left[\left(m-j+1\right)^{\alpha} -\left(m-j\right)^{\alpha}\right] , \quad \quad \quad \quad \quad \quad \quad \quad \quad \quad \quad 
\end{align}
\begin{align}
& \int^{t_{j+1}}_{t_{j}}\left(s-t_{j}\right) \left( t_{m+1}-s\right)^{\alpha-1}ds
=\frac{\left(\Delta t\right)^{\alpha+1}}{\alpha \left(\alpha+1\right)} \Bigg[\left(m-j+1\right)^{\alpha+1}
- \left(m-j\right)^{\alpha+1}-\left(\alpha+1\right)\left(m-j\right)^{\alpha}\Bigg],
\end{align}
\begin{multline}
\int^{t_{j+1}}_{t_{j}}\left(s-t_{j}\right)\left(s-t_{j+\frac{1}{2}}\right) \left( t_{m+1}-s\right)^{\alpha-1}ds=\frac{\left(\Delta t\right)^{\alpha+2}}{\alpha \left(\alpha+1\right)\left(\alpha+2\right)} \qquad \qquad \qquad \qquad \qquad \qquad \quad ,\\
\qquad \quad \times \Bigg [  \left(m-j+1\right)^{\alpha}\Bigg(2\left(m-j\right)^{2}-\frac{1}{2} \bigg[\left(\alpha-6\right)\left(m-j\right)+\left(\alpha-2\right)\bigg] \Bigg)\\
-\left(m-j\right)^{\alpha}\Bigg(2\left(m-j\right)^{2}+\frac{1}{2} \bigg[\left(3\alpha+6\right)\left(m-j\right)+\left(\alpha^{2}+3\alpha+2\right)\bigg] \Bigg)\Bigg], \qquad \,\, \,\,\,\,\,
\end{multline}
respectively.
Replacing them into  (\ref{eq3}),  we have the following 
\begin{multline}\label{eq4}
y_{m+1}=y_{0}+ \frac{1-\alpha}{AB(\alpha)}f\left(t_{m+1},y_{m+1}\right)\\
+\frac{\alpha \left(\Delta t\right)^{\alpha}}{AB(\alpha)\Gamma(\alpha+1)} \sum^{m}_{j=0}f\left(t_{j},y_{j}\right)  \left[\left(m-j+1\right)^{\alpha} -\left(m-j\right)^{\alpha}\right] \qquad \qquad \qquad \qquad \qquad \qquad \quad  \\ 
+\frac{2 \alpha \left(\Delta t\right)^{\alpha}}{AB(\alpha)\Gamma(\alpha+2)} \sum^{m}_{j=0}\Bigg(f\left(t_{j+\frac{1}{2}},y_{j+\frac{1}{2}}\right)-f\left(t_{j},y_{j}\right)\Bigg) \qquad  \qquad  \qquad \qquad  \qquad \qquad  \qquad \,\,\,\,\,\, \\
\times \Bigg[\left(m-j+1\right)^{\alpha+1}
- \left(m-j\right)^{\alpha+1}-\left(\alpha+1\right)\left(m-j\right)^{\alpha}\Bigg] \quad \,\,\,\, \qquad \qquad  \qquad \qquad \qquad  \qquad \, \, \\ 
+\frac{2\alpha  \left(\Delta t\right)^{\alpha}}{AB(\alpha)\Gamma(\alpha+3)} \sum^{m}_{j=0} \Bigg( f\left(t_{j},y_{j}\right)-2f\left(t_{j+\frac{1}{2}},y_{j+\frac{1}{2}}\right)+f\left(t_{j+1},y_{j+1}\right) \Bigg)\qquad   \qquad \qquad \quad \, \\
\times \Bigg [  \left(m-j+1\right)^{\alpha}\Bigg(2\left(m-j\right)^{2}-\frac{1}{2} \bigg[\left(\alpha-6\right)\left(m-j\right)+\left(\alpha-2\right)\bigg] \Bigg) \qquad   \qquad \qquad \quad \qquad  \,\,\, \, \\
-\left(m-j\right)^{\alpha}\Bigg(2\left(m-j\right)^{2}+\frac{1}{2} \bigg[\left(3\alpha+6\right)\left(m-j\right)+\left(\alpha^{2}+3\alpha+2\right)\bigg] \Bigg)\Bigg]. \qquad \qquad \qquad  \qquad \,\,\,\,\,\,\,\,\,
\end{multline}
In order to simplify the formulas to come,  let us define the expression
\begin{multline} 
\Upsilon_{p}= 
\frac{ \left(\Delta t\right)^{\alpha}}{\Gamma(\alpha+1)} \sum^{p}_{j=0}f\left(t_{j},y_{j}\right)  \left[\left(m-j+1\right)^{\alpha} -\left(m-j\right)^{\alpha}\right] \qquad \qquad \qquad \qquad \qquad \qquad \quad  \\ 
+\frac{2  \left(\Delta t\right)^{\alpha}}{\Gamma(\alpha+2)} \sum^{p}_{j=0}\left(f\left(t_{j+\frac{1}{2}},y_{j+\frac{1}{2}}\right)-f\left(t_{j},y_{j}\right)\right)
\Bigg[\left(m-j+1\right)^{\alpha+1}- \left(m-j\right)^{\alpha+1}-\left(\alpha+1\right)\left(m-j\right)^{\alpha}\Bigg]  \\ 
+\frac{ \left(\Delta t\right)^{\alpha}}{\Gamma(\alpha+3)} \sum^{p}_{j=0} \Bigg( 2f\left(t_{j},y_{j}\right)-4f\left(t_{j+\frac{1}{2}},y_{j+\frac{1}{2}}\right)+f\left(t_{j+1},y_{j+1}\right) \Bigg)\qquad   \qquad \qquad \quad \qquad \qquad \qquad  \quad \,\,\, \\
\times \Bigg [  \left(m-j+1\right)^{\alpha}\Bigg(2\left(m-j\right)^{2}-\frac{1}{2} \bigg[\left(\alpha-6\right)\left(m-j\right)+\left(\alpha-2\right)\bigg] \Bigg) \qquad   \qquad \qquad \quad \qquad  \qquad  \qquad  \,\,\, \,\,\,\, \\
-\left(m-j\right)^{\alpha}\Bigg(2\left(m-j\right)^{2}+\frac{1}{2} \bigg[\left(3\alpha+6\right)\left(m-j\right)+\left(\alpha^{2}+3\alpha+2\right)\bigg] \Bigg)\Bigg], \qquad \qquad \qquad  \quad \quad \qquad  \,\,\,\,\,\,\,\,\,\,
\end{multline}
with the convention 
\begin{align}
&\Upsilon_{0}=0.
\end{align}
Using the above  notation, (\ref{eq4}) can be rewritten in the form 
\begin{multline}\label{eqcorrector}
y_{m+1}=y_{0}+ \frac{1-\alpha}{AB(\alpha)}f\left(t_{m+1},y_{m+1}\right) +\frac{\alpha }{AB(\alpha)} \Upsilon_{m-1} \\
+\frac{\alpha \left(\Delta t\right)^{\alpha}}{AB(\alpha)\Gamma(\alpha+1)} f\left(t_{m},y_{m}\right) \qquad \qquad \qquad \qquad \qquad \qquad \quad \qquad \qquad \qquad \qquad \qquad \qquad \quad  \\
+\frac{2 \alpha \left(\Delta t\right)^{\alpha}}{AB(\alpha)\Gamma(\alpha+2)}  \Bigg( f\left( t_{m+\frac{1}{2}}, y_{m+\frac{1}{2}}\right)-f\left(t_{m},y_{m}\right) \Bigg) \qquad \qquad \qquad \qquad \qquad \qquad \quad \qquad  \,\,  \,\, \, \\
-\frac{ \alpha \left(\alpha-2\right) \left(\Delta t\right)^{\alpha}}{2 AB(\alpha)\Gamma(\alpha+3)} \Bigg( 2f\left(t_{m},y_{m}\right)-4f\left(t_{m+\frac{1}{2}},y_{m+\frac{1}{2}}\right)+f\left(t_{m+1},y_{m+1}\right) \Bigg). \qquad \qquad \qquad \quad \,\,\, \,\, \,
\end{multline}
Formula (\ref{eqcorrector}) will serve as our implicit part,  i.e. the corrector.
To derive an explicit formula for approximating $y_{j+\frac{1}{2}}$
\begin{multline}\label{eq6}
y_{j+\frac{1}{2}}= y_{0}+ \frac{1-\alpha}{AB(\alpha)}f\left(t_{j},y_{j}\right)\\
+\frac{\alpha}{AB(\alpha)\Gamma(\alpha)}  \Bigg\{ \sum^{j-1}_{k=0} \Bigg[  \int^{t_{k+1}}_{t_{k}} \hat{f_{k}}\left(s,y\left(s\right)\right)\left(t_{j+\frac{1}{2}}-s\right)^{\alpha-1}ds \Bigg] \qquad \qquad \qquad \qquad \qquad \qquad \,\,\,\,\,   \\
+ \int^{t_{j+\frac{1}{2}}}_{t_{j}} \bar{f_{k}}\left(s,y\left(s\right)\right)\left(t_{j+\frac{1}{2}}-s\right)^{\alpha-1}ds\Bigg\},\qquad \qquad\qquad \qquad\qquad \qquad \qquad \qquad \qquad \,\qquad \quad \,\,\,\,\,\,
\end{multline}
where $\hat{f_{k}}\left(s,y\left(s\right)\right)$ and $\bar{f_{k}}\left(s,y\left(s\right)\right)$ is the piecewise linear Newton interpolation given by 
\begin{equation}
\begin{cases}
&\hat{f_{k}}\left(s,y\left(s\right)\right)= f\left(t_{k},y_{k}\right)+\frac{ f\left(t_{k+1},y_{k+1}\right)-f\left(t_{k},y_{k}\right)}{\Delta t} \left(s-t_{k}\right) ,   \qquad   \qquad    \qquad     \\
& \bar{f_{k}}\left(s,y\left(s\right)\right)=  f\left(t_{j},y_{j}\right).
\end{cases}
\end{equation}
Thus if we put this polynomial in  formula (\ref{eq6}),  we have the following 
\begin{multline}\label{eq7}
y_{j+\frac{1}{2}}=y_{0}+ \frac{1-\alpha}{AB(\alpha)}f\left(t_{j},y_{j}\right)\\
+\frac{\alpha}{AB(\alpha)\Gamma(\alpha)} \sum^{j-1}_{k=0} f\left(t_{k},y_{k}\right) \int^{t_{k+1}}_{t_{k}}\left(t_{j+\frac{1}{2}}-s\right)^{\alpha-1}ds \qquad \qquad \qquad \qquad \qquad \qquad \quad \qquad \,\,\,\,   \\ 
+\frac{ \alpha}{AB(\alpha)\Gamma(\alpha)} \sum^{j-1}_{k=0}\frac{f\left(t_{k+1},y_{k+1}\right)-f\left(t_{k},y_{k}\right)}{\Delta t} \int^{t_{k+1}}_{t_{k}}\left(s-t_{k}\right)\left(t_{j+\frac{1}{2}}-s\right)^{\alpha-1}ds \qquad \qquad \,\,\,\,\,\,    \\ 
+\frac{ \alpha}{AB(\alpha)\Gamma(\alpha)} f\left(t_{j},y_{j}\right)\int^{t_{j+\frac{1}{2}}}_{t_{j}} \left(t_{j+\frac{1}{2}}-s\right)^{\alpha-1}ds.\qquad \qquad \qquad \qquad \qquad \qquad \quad \qquad \quad \qquad \,\, \,\,\,
\end{multline}
When calculating the integral in (\ref{eq7})
\begin{multline}
\quad \quad  \,\,  \,\,    \int^{t_{k+1}}_{t_{k}}\left(t_{j+\frac{1}{2}}-s\right)^{\alpha-1}ds=\frac{\left(\Delta t\right)^{\alpha}}{\alpha} \bigg[\left(j-k+\frac{1}{2}\right)^{\alpha}- \left(j-k-\frac{1}{2}\right)^{\alpha}\bigg] , \\
\int^{t_{k+1}}_{t_{k}}\left(s-t_{k}\right)\left(t_{j+\frac{1}{2}}-s\right)^{\alpha-1}ds = \frac{\left(\Delta t\right)^{\alpha+1}}{\alpha \left(\alpha+1\right)} \times \qquad \quad \quad \quad \qquad \qquad \quad \quad \quad \qquad \qquad \quad \quad \quad \qquad \qquad \quad \quad \quad  \\  
\qquad \qquad \quad \quad \quad \qquad  \quad \quad \quad   \Bigg[\left(j-k+\frac{1}{2}\right)^{\alpha+1}
- \left(j-k-\frac{1}{2}\right)^{\alpha+1}-\left(\alpha+1\right)\left(j-k-\frac{1}{2}\right)^{\alpha}\Bigg],\\
\qquad \, \,\, \, \int^{t_{j+\frac{1}{2}}}_{t_{j}} \left(t_{j+\frac{1}{2}}-s\right)^{\alpha-1}ds= \frac{\left(\frac{\Delta t}{2}\right)^{\alpha}}{\alpha}, \qquad \qquad \qquad \qquad \qquad \qquad \quad \qquad \quad \qquad \,\, \qquad \qquad \quad \quad \quad  \,\,  \,\,  
\end{multline}
and putting these calculation into (\ref{eq7}) we can obtain the following 
\begin{multline}
y_{j+\frac{1}{2}}=y_{0}+ \frac{1-\alpha}{AB(\alpha)}f\left(t_{j},y_{j}\right)\\
+\frac{\alpha}{AB(\alpha)\Gamma(\alpha)} \sum^{j-1}_{k=0} f\left(t_{k},y_{k}\right)\Bigg\{ \frac{\left(\Delta t\right)^{\alpha}}{\alpha} \bigg[\left(j-k+\frac{1}{2}\right)^{\alpha}- \left(j-k-\frac{1}{2}\right)^{\alpha}\bigg]\Bigg\}  \quad \qquad \qquad  \qquad   \\ 
+\frac{ \alpha}{AB(\alpha)\Gamma(\alpha)} \sum^{j-1}_{k=0}\frac{f\left(t_{k+1},y_{k+1}\right)-f\left(t_{k},y_{k}\right)}{\Delta t} \quad \qquad \qquad  \qquad \quad \qquad  \qquad \qquad  \qquad \qquad \quad \,\,\\
\times \Bigg\{ \frac{\left(\Delta t\right)^{\alpha+1}}{\alpha \left(\alpha+1\right)}\Bigg[\left(j-k+\frac{1}{2}\right)^{\alpha+1}
- \left(j-k-\frac{1}{2}\right)^{\alpha+1}-\left(\alpha+1\right)\left(j-k-\frac{1}{2}\right)^{\alpha}\Bigg]\Bigg\}\qquad \qquad \\ 
+\frac{ \alpha}{AB(\alpha)\Gamma(\alpha)} f\left(t_{j},y_{j}\right)\Bigg\{ \frac{\left(\frac{\Delta t}{2}\right)^{\alpha}}{\alpha} \Bigg\}.\qquad \qquad \qquad \qquad \qquad \qquad \quad \qquad \quad \qquad \qquad \quad \qquad \quad  \,\,\,\,\,\,
\end{multline}
Simplifying and rearranging the terms leads to
\begin{multline}\label{eqmid}
y_{j+\frac{1}{2}}=y_{0}+ \frac{1-\alpha}{AB(\alpha)}f\left(t_{j},y_{j}\right)\\
+\frac{\alpha \left(\Delta t\right)^{\alpha}}{AB(\alpha)\Gamma(\alpha+1)} \sum^{j-1}_{k=0} f\left(t_{k},y_{k}\right) \bigg[\left(j-k+\frac{1}{2}\right)^{\alpha}- \left(j-k-\frac{1}{2}\right)^{\alpha}\bigg] \quad \qquad \qquad  \qquad \quad \quad \,\,\,\,\\ 
+\frac{ \alpha \left(\Delta t\right)^{\alpha} }{AB(\alpha)\Gamma(\alpha+2)} \sum^{j-1}_{k=0}(f\left(t_{k+1},y_{k+1}\right)-f\left(t_{k},y_{k}\right)) \quad \qquad \qquad  \qquad \quad \qquad  \qquad \qquad  \qquad \quad  \,\,\\
\times\Bigg[\left(j-k+\frac{1}{2}\right)^{\alpha+1}
- \left(j-k-\frac{1}{2}\right)^{\alpha+1}-\left(\alpha+1\right)\left(j-k-\frac{1}{2}\right)^{\alpha}\Bigg]\qquad \qquad \qquad \qquad \quad\,\, \,\\ 
+\frac{ \alpha \left(\frac{\Delta t}{2}\right)^{\alpha}}{AB(\alpha)\Gamma(\alpha+1)} f\left(t_{j},y_{j}\right).\qquad \qquad \qquad \qquad \qquad \qquad \quad \qquad \quad \qquad \qquad \quad  \quad \quad \qquad \qquad  \,\,\,\,\,
\end{multline}
To obtain our predictor formula, we consider the following non-linear fractional differential equation
\begin{multline}\label{eq8}
y^P_{m+1}=y_{0}+ \frac{1-\alpha}{AB(\alpha)}f\left(t_{m+\frac{1}{2}},y_{m+\frac{1}{2}}\right)+\frac{\alpha}{AB(\alpha)\Gamma(\alpha)} \Bigg[  \sum^{m}_{j=0} \int^{t_{j+1}}_{t_{j}}\bar{\bar f}_{j}\left(s,y\left(s\right)\right)\left(t_{m+1}-s\right)^{\alpha-1}ds \Bigg] ,\quad
\end{multline}
where
\begin{equation}\label{eq9}
\bar{\bar f}_{j}\left(s,y\left(s\right)\right)=\left \{ 
\begin{array}{lll}
f\left(t_{j},y_{j}\right)+\frac{2 \left[f\left(t_{j+\frac{1}{2}},y_{j+\frac{1}{2}}\right)-f\left(t_{j},y_{j}\right)\right]}{\Delta t} \left(s-t_{j}\right) & \text{if} & s\in \left[t_{j},t_{j+\frac{1}{2}} \right]
, \\
f\left(t_{j+\frac{1}{2}},y_{j+\frac{1}{2}}\right) & \text{if} & s\in \left[t_{j+\frac{1}{2}},t_{j+1} \right].\mathbb{\medskip }
\end{array}%
\right. 
\end{equation}
By using the Newton polynomial (\ref{eq9}),  formula  (\ref{eq8})  becomes
\begin{eqnarray}\label{eq10}
y^{P}_{m+1}&=&y_{0}+ \frac{1-\alpha}{AB(\alpha)}f\left(t_{m+\frac{1}{2}},y_{m+\frac{1}{2}}\right) \notag\\
&&+\frac{\alpha}{AB(\alpha)\Gamma(\alpha)} \sum^{m}_{j=0} f\left(t_{j},y_{j}\right) \int^{t_{j+\frac{1}{2}}}_{t_{j}}\left(t_{m+1}-s\right)^{\alpha-1}ds  \notag \\
&&+\frac{2 \alpha}{AB(\alpha)\Gamma(\alpha)} \sum^{m}_{j=0}\left( \frac{f\left(t_{j+\frac{1}{2}},y_{j+\frac{1}{2}}\right)-f\left(t_{j},y_{j}\right)}{\Delta t}\right)  \int^{t_{j+\frac{1}{2}}}_{t_{j}}\left(s-t_{j}\right)\left(t_{m+1}-s\right)^{\alpha-1}ds \notag\\
&&+\frac{\alpha}{AB(\alpha)\Gamma(\alpha)} \sum^{m}_{j=0} f\left(t_{j+\frac{1}{2}},y_{j+\frac{1}{2}}\right) \int_{t_{j+\frac{1}{2}}}^{t_{j+1}}\left(t_{m+1}-s\right)^{\alpha-1}ds.
\end{eqnarray}
The integrals in (\ref{eq10}) can be calculated as
\begin{equation}
\int^{t_{j+\frac{1}{2}}}_{t_{j}}\left(t_{m+1}-s\right)^{\alpha-1}ds=\frac{(\Delta t)^\alpha}{\alpha}\left[\left( m-j+1\right) ^\alpha -\left( m-j+\frac{1}{2}\right) ^\alpha \right] ,  \quad \quad \qquad \qquad  
\end{equation}
\begin{eqnarray}
\int^{t_{j+\frac{1}{2}}}_{t_{j}}\left(s-t_{j}\right)\left(t_{m+1}-s\right)^{\alpha-1}ds&=&\frac{(\Delta t)^{\alpha+1}}{\alpha(\alpha+1)} \left(m+\alpha j+1 \right)\left(m-j+1 \right)^\alpha \notag\\
&&-\frac{(\Delta t)^{\alpha+1}}{\alpha(\alpha+1)}\left( m+\alpha j+\frac{\alpha+2}{2}\right)\left( m-j+\frac{1}{2}\right)^\alpha \notag\\
&&+ \frac{(\Delta t)^{\alpha+1}}{\alpha(\alpha+1)}\left[(\alpha+1)j\left(m-j+\frac{1}{2} \right)^\alpha - (\alpha+1)j\left(m-j+1 \right)^\alpha \right] ,
\end{eqnarray}
\begin{equation}
\int_{t_{j+\frac{1}{2}}}^{t_{j+1}}\left(t_{m+1}-s\right)^{\alpha-1}ds=\frac{(\Delta t)^{\alpha}}{\alpha}\left[\left(m-j+\frac{1}{2}\right)^\alpha -\left(m-j\right)^\alpha \right] .  \quad \quad  \quad \quad \qquad \quad  \,\,\,
\end{equation}
By substituting these calculations into  (\ref{eq9}),  we
obtain%
\begin{eqnarray}\label{eqpredictor}
y^{P}_{m+1}&=&y_{0}+ \frac{1-\alpha}{AB(\alpha)}f\left(t_{m+\frac{1}{2}},y_{m+\frac{1}{2}}\right) \notag\\
&&+\frac{(\Delta t)^\alpha}{AB(\alpha)\Gamma(\alpha)} \sum^{m}_{j=0} f\left(t_{j},y_{j}\right) \left[\left( m-j+1\right) ^\alpha -\left( m-j+\frac{1}{2}\right) ^\alpha \right]   \notag \\
&&+\frac{(\Delta t)^{\alpha}}{AB(\alpha)\Gamma(\alpha)} \sum^{m}_{j=0} f\left(t_{j+\frac{1}{2}},y_{j+\frac{1}{2}}\right) \left[\left(m-j+\frac{1}{2}\right)^\alpha -\left(m-j\right)^\alpha \right]  \\
&&+\frac{2 \alpha(\Delta t)^{\alpha}}{AB(\alpha)\Gamma(\alpha+2)} \sum^{m}_{j=0}\left(f\left(t_{j+\frac{1}{2}},y_{j+\frac{1}{2}}\right)-f\left(t_{j},y_{j}\right)\right)  \mathcal{C}(\alpha, j, m) ,\notag
\end{eqnarray}
where
\begin{eqnarray}
\mathcal{C}(\alpha, j, m)&=& \left(m+\alpha j+1 \right)\left(m-j+1 \right)^\alpha -\left(m+\alpha j+\frac{\alpha+2}{2} \right) \left(m-j+\frac{1}{2} \right)^\alpha \notag \\
&&+\left(\alpha+1 \right) j\left(m-j+\frac{1}{2} \right)^\alpha -\left(\alpha+1 \right)j \left(m-j+1 \right)^\alpha. \notag 
\end{eqnarray}
Formula (\ref{eqpredictor}) will serve as our explicit, i.e the predictor.  In each iteration,  the predictor in (\ref{eqpredictor})is calculated and then corrected by means of the corrector formula (\ref{eqcorrector}) .

\section{Numerical Illustrations and Simulation}
In this section,  we will gives some simulation results for the fractional  predictor-corrector numerical scheme (\ref{eqcorrector}) -(\ref{eqmid})  and (\ref{eqpredictor}), the predictor-corrector numerical method introduced by Baleanu-Jajarmi \cite{
baleanu2018nonlinear} and the numerical scheme introduced by Toufik-Atangana \cite{
Toufik2017} to show the efficiency  and accuracy of our new method . The experimental order of convergence is computed by 
\begin{equation}
EOC= log_{2} \frac{E\left(\frac{N}{2}\right)}{E\left(N\right)},
\end{equation}
where $E\left(N\right)=$max$_{1\leq i\leq N}  |y\left(t_{i}\right)-y_{i}|$ .

\textbf{Example 1} Consider the following fractional Initial value problem 
\begin{equation}\label{Exm1}
\begin{cases}
&{}^{ABC}_{\quad \,\, 0}D^{\alpha}_{t} y\left(t\right)= t^{n},   \qquad   \qquad    \qquad       0< t \leq 1, \\ \\
& \qquad  \quad \, y\left(0\right)= 1,    \qquad    \qquad   \qquad      \,\,\,       
\end{cases}
\end{equation}
where $\alpha\in(0,1)$ and $n$ is a positive integer.  According  \cite{baleanu2018nonlinear},  the exact solution of ($32$) is given by 
\begin{align}
&y\left(t\right)=y_{0}+\frac{1-\alpha}{B\left(\alpha\right)}t^{n}+\frac{\alpha \Gamma\left(n+1\right)}{B\left(\alpha\right)\Gamma\left(\alpha+n+1\right)}t^{\alpha+n}.
\end{align}
Example $1$ is numerically solved by present method  (\ref{eqcorrector}) -(\ref{eqmid})  and (\ref{eqpredictor}) and other method for various values of $\alpha$.  
Table $1$ and $2$ show the comparison of absolute error of different numerical scheme for example $1$ with 
$t \in \left[0,1\right]$ and various value of parameter $\alpha$. The results in table $2$ show that the orders of convergence for  Baleanu-Jajarmi numerical scheme are roughly $2$,  and for present method are roughly $3.5$ particularly when steps size h is small($h=\frac{1}{N}$).  
From figure $1$,  we note that the present numerical scheme  (\ref{eqcorrector}) -(\ref{eqmid})  and (\ref{eqpredictor}) agree well with the exact solution . 
From these tables,  we see that the proposed numerical scheme (\ref{eqcorrector}) -(\ref{eqmid})  and (\ref{eqpredictor}) achieve a lower error than the Toufik-Atangana and Baleanu-Jajarmi numerical scheme .  Also,  we observe that the solution obtained by proposed scheme increases rapidly to the exact solution following the increase in  $\alpha$ . 
Those figure and tables indicate  the efficacy of the present predictor-corrector numerical method
 (\ref{eqcorrector}) -(\ref{eqmid})  and (\ref{eqpredictor}) for solving IVP.
\begin{table}
	\caption{ The absoluate error of various numerical methods for problem $(32)$ with $n=2$, N=$40$ and $t \in \left[0,1\right]$ using various values of $\alpha$.}
	\begin{center}
	\begin{tabular}{lllllllll}
		\hline\noalign{\smallskip}
		Methods   & \qquad $\alpha=0.5$ & \qquad $\alpha=0.7$ & \qquad $\alpha=0.99$ & \qquad $\alpha=1$\\
		\noalign{\smallskip}
		\hline
		\noalign{\smallskip}
The proposed numerical scheme  (\ref{eqcorrector}) -(\ref{eqmid})  and (\ref{eqpredictor})
			& 
			\qquad $2.0 e{-16}$
			& 
			\qquad $4.0 e{-16}$
			& 
			\qquad $4.0 e{-16}$
			& 
			\qquad $2.0 e{-16}$
			\parbox[t]{1cm}{
			\raggedright}\\
 Baleanu-Jajarmi numerical scheme \cite{baleanu2018nonlinear}
			& 
			\qquad $5.6 e^{-5}$
			& 
			\qquad $7.9 e{-5}$
			& 
			\qquad $4.7 e{-5}$
			& 
			\qquad $1.0 e{-4}$
			\parbox[t]{1cm}{
			\raggedright}\\
Toufik-Atangana numerical scheme \cite{Toufik2017} 
			& 
			\qquad $2.5 e^{-2}$
			& 
			\qquad $0.02$
			& 
			\qquad $5.4 e^{-3}$
			& 
			\qquad $5.0 e^{-4}$
			\parbox[t]{1cm}{
			\raggedright}\\						
		\hline
	\end{tabular}
\end{center}
\end{table}

\begin{table}[h]
        \caption{The absolute error (AE), experimental order of convergence (EOC) and CPU time in seconds (CTs) of the numerical scheme  (\ref{eqcorrector})-(\ref{eqmid})  and (\ref{eqpredictor})   and the numerical scheme Baleanu-Jajarmi  \cite{baleanu2018nonlinear} for example 1 with $n=3$ and $t \in \left[0,1\right]$.	 }
\begin{tabular}{l ccccc ccccc  }
		\hline
\multicolumn{1}{c}{Method} & \multicolumn{1}{c}{N}  & \multicolumn{3}{c}{$\alpha=0.5$} & \multicolumn{3}{c}{$\alpha=0.7$} & \multicolumn{3}{c}{$\alpha=0.9$} \\
&& AE  & EOC  & CTs  & AE  & EOC  & CTs   & AE  & EOC  & CTs   \\
		\hline
Proposed numerical 					&10     & $2.1e{-6}$ & $-$ & $0.08$  & 		 $1.3e{-6}$  & $-$ & $6.7e{-3}$& $4.2e{-7}$  & $-$ & $1.3e{-3}$   \\
method   (\ref{eqcorrector})-(\ref{eqmid})  and (\ref{eqpredictor})        &20    & $2.0e{-7}$ & $3.45$ & $0.01$   & $1.1e{-7}$ & $3.61$ & $5.7e{-3}$& $3.1e{-8}$  & $3.74$ & $7.6e{-3}$  \\
		  										&40    & $1.8e{-8}$& $3.46$ & $0.04$   & $9.1e{-9}$ & $3.63$& $0.036$ & $2.3e{-9}$  & $3.77$ & $0.039$ \\
												&80    & $1.6e{-9}$& $3.74$ & $0.19$   & $7.2e{-10}$  &$3.64$ & $0.26$ & $1.6e{-10}$  & $3.79$ & $0.25$\\
		\hline
Baleanu-Jajarmi 														  &10    		 & $1.6e{-3}$ & $-$ & $0.02$   & $2.1e{-3}$ & $-$ & $1.3e{-3}$ & $2.4e{-3}$ & $-$ & $7.0e{-4}$   \\
numerical	method  \cite{baleanu2018nonlinear}		       &20        & $4.3e{-4}$ & $1.95$ & $2.1e{-3}$ & $5.5e{-4}$ & $1.98$ & $1.1e{-3}$    & $6.1e{-4}$ & $1.99$ & $1.7e{-3}$   \\
		  										                                     &40    		& $1.1e{-4}$ & $1.96$ & $8.0e{-3}$ &  $1.3e{-4}$ & $1.98$ & $2.9e{-3}$ & $1.5e{-4}$ & $1.99$ & $2.8e{-3}$   \\
																				 &80    		& $2.8e{-5}$ & $1.97$ & $8.2e{-3}$  & $3.5e{-5}$ & $1.99$ & $8.3e{-3}$  & $3.8e{-5}$ & $1.99$ & $8.0e{-3}$   \\
		\hline
\end{tabular}
\end{table}

\begin{center}
	\begin{figure}[!h]
		\includegraphics[width=17cm]{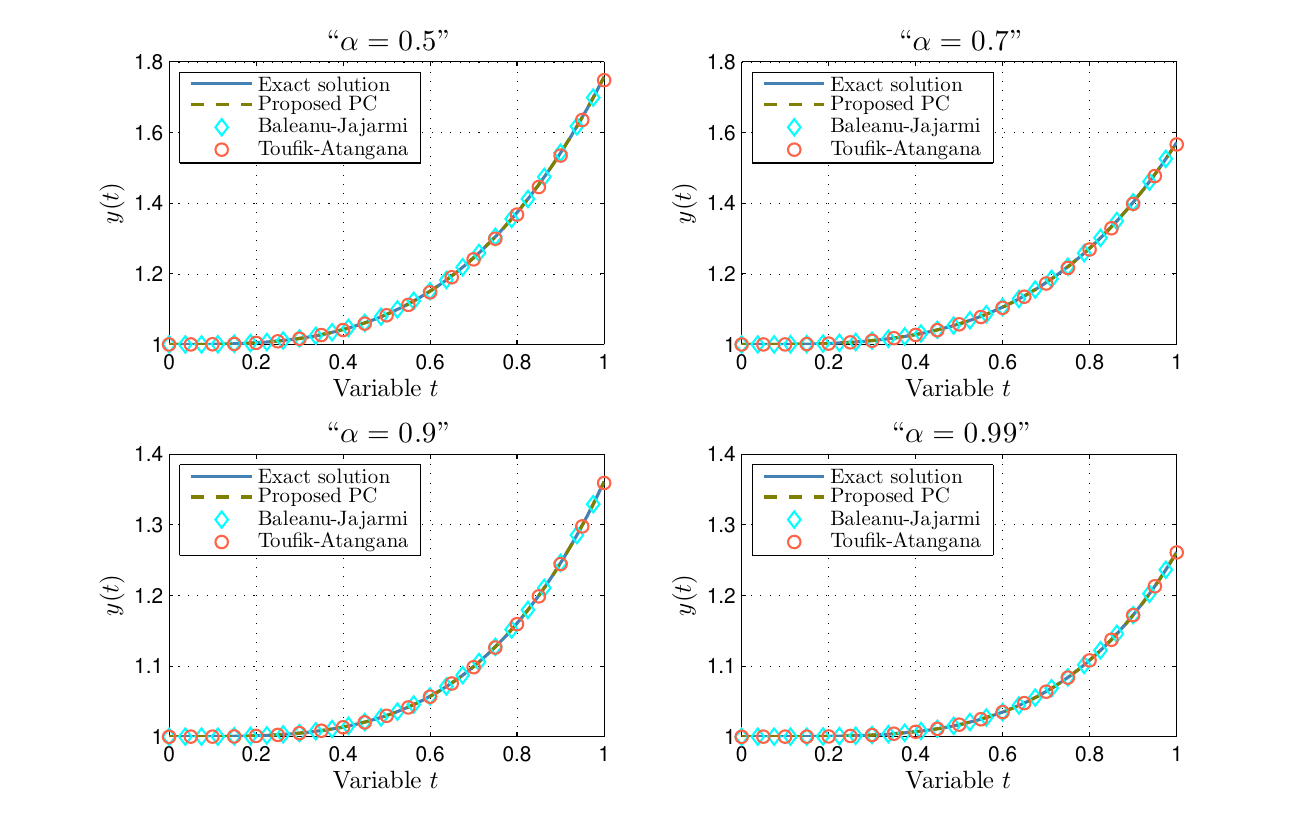}
		\caption{Comparison of the exact and the numerical solutions of various method for (\ref{Exm1}), when $\alpha\in{0.5,0.7,0.9,0.99}$, $n=3$, $N=160$ and $t \in \left[0,1\right]$. }
		\label{simulation_1}
	\end{figure}
\end{center}

\textbf{Example 2} Consider the following fractional differential equation
\begin{equation}
\begin{cases}
&{}^{ABC}_{\quad \,\, 0}D^{\alpha}_{t} y\left(t\right)+ y\left(t\right)= t,   \qquad   \qquad    \qquad       0< t \leq 1, \\ \\
& \qquad  \quad \, y\left(0\right)= 0,   \qquad    \qquad   \qquad     \qquad    \quad     \
\end{cases}
\end{equation}
where $\alpha \in (0,1)$. According to  \cite{baleanu2018nonlinear}, the exact solution of ($34$) is given by 
\begin{align}\label{Exm2}
&y\left(t\right)=\frac{1}{B\left(\alpha\right)+1-\alpha} \Bigg\{ \left(1-\alpha\right)t   \textbf{E}_{\alpha,2}\left( -\frac{\alpha}{B\left(\alpha\right)+1-\alpha} t^{\alpha}\right) +\alpha t^{\alpha+1}    \textbf{E}_{\alpha,\alpha+2}\left(- \frac{\alpha}{B\left(\alpha\right)+1-\alpha} t^{\alpha}\right)  \Bigg\}.
\end{align}
In this example,  we also apply proposed predictor-corrector  method (\ref{eqcorrector}) -(\ref{eqmid})  and (\ref{eqpredictor})and other method at time $t = 1$ with different steps and different values of $\alpha$ to get the numerical solutions,  the results are obtained shown in Tables $3$.  Table $3$ also show that  the orders of convergence for  present method are smaller than  Baleanu-Jajarmi numerical scheme.  In Fig. 2, where $\alpha=0.25,0.3,0.5,0.7$ and $N = 160$,  the numerical solution determined by  (\ref{eqcorrector}) -(\ref{eqmid})  and (\ref{eqpredictor}) and the exact solution are compared very well.  Once again, the present scheme has better results than other methods.
\begin{table}[h]
\caption{The absolute error (AE), experimental order of convergence (EOC) and CPU time in seconds (CTs) of the numerical scheme  (\ref{eqcorrector})-(\ref{eqmid})  and(\ref{eqpredictor}) and the numerical scheme Baleanu-Jajarmi  \cite{baleanu2018nonlinear} for example 2.}
\begin{tabular}{l ccccc ccccc  }
		\hline
\multicolumn{1}{c}{Method} & \multicolumn{1}{c}{N}  & \multicolumn{3}{c}{$\alpha=0.3$} & \multicolumn{3}{c}{$\alpha=0.5$} & \multicolumn{3}{c}{$\alpha=0.7$} \\
&& AE  & EOC  & CTs  & AE  & EOC  & CTs   & AE  & EOC  & CTs   \\
		\hline
Proposed numerical 				                                    	  &10     & $9.4 e{-3}$ & $-$ & $7.6 e{-3}$  & $9.4 e{-3}$  & $-$ & $3.0 e{-4}$                     & $4.8 e{-3}$  & $-$ & $5.0 e{-4}$   \\
method   (\ref{eqcorrector})-(\ref{eqmid})  and (\ref{eqpredictor})     									  &20    & $5.0 e{-3}$ & $0.89$ & $4.8 e{-3}$   & $4.6 e{-3}$ & $1.02$ & $1.3 e{-3}$                & $2.2 e{-3}$  & $1.12$ & $2.5 e{-3}$  \\
		  																		 &40    & $2.7 e{-3}$& $0.92$ & $0.03$   & $2.3 e{-3}$ & $1.02$& $9.3 e{-3}$                     & $1.0 e{-3}$  & $1.08$ & $1.7 e{-2}$ \\
																				&80    & $1.4 e{-3}$& $0.94$ & $0.13$   & $1.1 e{-3}$  &$1.02$& $7.0 e{-2}$                     & $5.0 e{-4}$  & $1.05$ & $0.14$\\
		\hline
Baleanu-Jajarmi 														  &10    		 & $0.11$ & $-$ & $4.8 e{-3}$   & $3.4 e{-2}$ & $-$ & $3.0 e{-4}$                        & $1.2 e{-2}$ & $-$ & $3.0 e{-4}$   \\
numerical	method  \cite{baleanu2018nonlinear}		       &20        & $0.05$ & $1.19$ & $1.4 e{-2}$ & $1.5 e{-2}$ & $1.17$ & $5.0 e{-4}$                    & $5.0 e{-3}$ & $1.22$ & $5.0e{-4}$   \\
		  										                                     &40    		& $0.02$ & $1.15$ & $2.4 e{-3}$ &  $6.9 e{-3}$ & $1.12$ & $1.3 e{-3}$                 & $2.0 e{-3}$ & $1.15$ & $1.5 e{-3}$   \\
																				 &80    		& $9.6 e{-3}$ & $1.11$ & $3.5 e{-3}$  & $3.2 e{-3}$ & $1.08$ & $2.8 e{-3}$            & $1.0 e{-3}$ & $1.10$ & $5.0 e{-3}$   \\
		\hline
\end{tabular}
\end{table}

\begin{center}
	\begin{figure}[!h]
		\includegraphics[width=17cm]{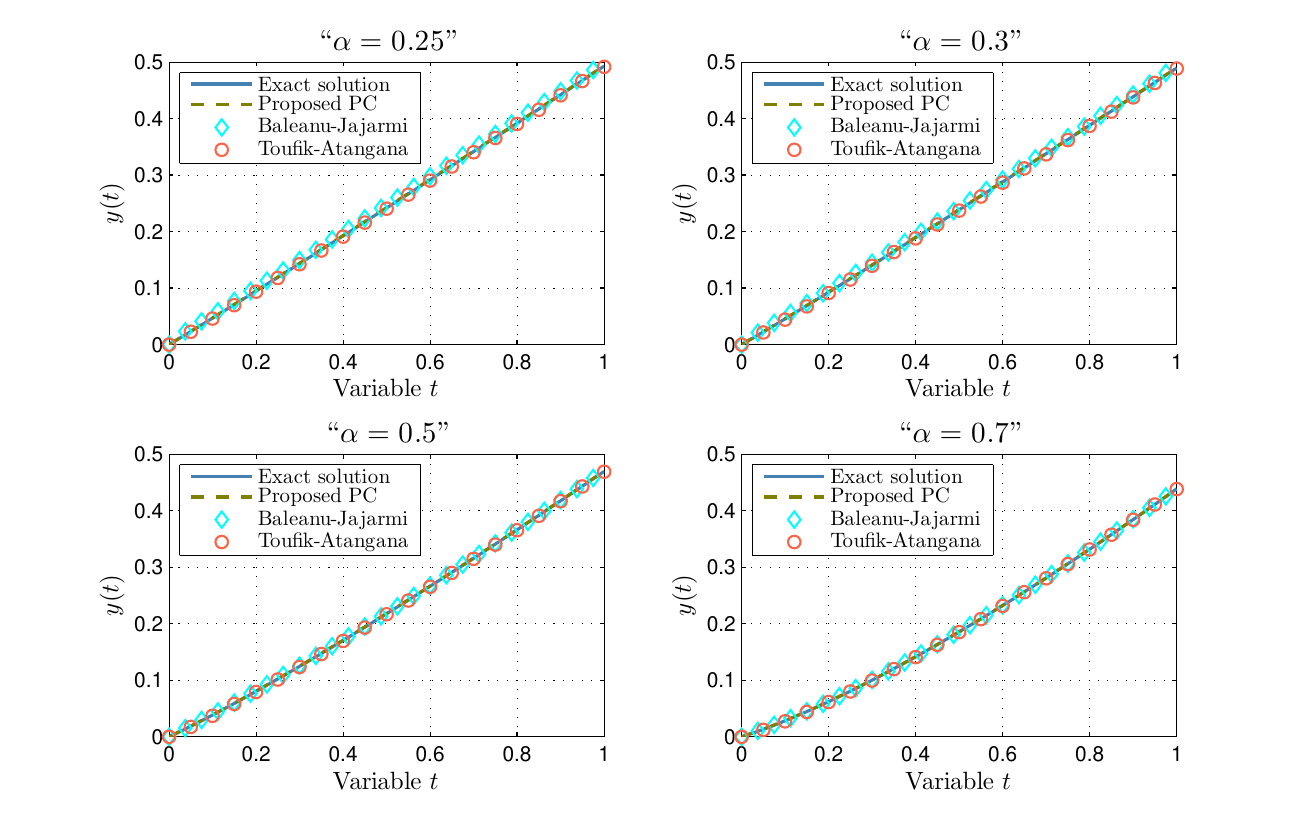}
		\caption{Comparison of the exact and the numerical solutions of various method for (\ref{Exm2}), when $\alpha\in{0.25,0.3,0.5,0.7}$, $n=3$, $N=160$ and $t \in \left[0,1\right]$. }
		\label{simulation_2}
	\end{figure}
\end{center}
\newpage
\section{Conclusion }
In this paper we propose  a new predictor-corrector method for solving fractional initial value problem based on Newton quadratic interpolation which includes an operator of the type of the Atangana–Baleanu fractional derivative.  The solutions obtained were compared with the exact solution,  which were in excellent agreement.  The numerical scheme effectively and accurately solves fractional initial value problems, as demonstrated by numerical examples.  Also,  based on numerical results,  the proposed method is shown to much more accurate than other methods.  Our results can be also extended to other definitions of fractional derivative.  The scope of numerical analysis of fractional initial value problem with a variety of fractional derivative definitions will be broadened by this method.

\end{document}